\documentclass[12pt, a4paper]{article}
\usepackage[latin1]{inputenc}
\usepackage[english]{babel}
\usepackage{indentfirst}
\usepackage{amsxtra}
\usepackage{amsmath}
\usepackage{amsthm}
\usepackage{amstext}
\usepackage{amsfonts}
\usepackage{textcomp}
\usepackage{amssymb}
\usepackage{amscd}
\usepackage[all]{xy}
\usepackage[dvips]{graphicx}
\usepackage{setspace}
\usepackage[all]{xy}

\newcommand {\demo}{\hskip -0.6cm{\bf Proof.  }}
\newcommand {\fim}{\hfill{$\square$}\vskip 1pc}

\newcommand {\R}{\mathbb{R}}
\newcommand {\N}{\mathbb{N}}

\newcommand {\E}{\mathbb{E}}

\newcommand{\af}[1]{\vspace{0.3cm}\hskip -0.6cm {\it Step #1: }\vspace{0.3cm}}

\newcommand{\fimaf}{\vspace{0.3cm}}

\newcommand{\funcao}[5]{\begin{array}{lrcl}
#1:&\!\!\!#2 & \rightarrow & #3 \\
  &\!\!\! #4 & \mapsto & #5
\end{array}}

\newtheorem{teorema}{Theorem}[section]

\newtheorem{corolario}[teorema]{Corollary}
\newtheorem{definicao}[teorema]{Definition}
\newtheorem{proposicao}[teorema]{Proposition}
\newtheorem{exemplo}[teorema]{Example}
\newtheorem{obs}[teorema]{Remark}

\begin{document}

\onehalfspace

\title{On the Representations of Leavitt path algebras}
\maketitle
\begin{center}
{\large Daniel Gonçalves and Danilo Royer}\\
\end{center}  
\vspace{8mm}

\begin{abstract}
Given a graph $E$ we define $E$-algebraic branching systems, show their existence and how they induce representations of the associated Leavitt path algebra. We also give sufficient conditions to guarantee faithfulness of the representations associated to E-algebraic branching systems and to guarantee equivalence of a given representation (or a restriction of it) to a representation arising from an E-algebraic branching system.
\end{abstract}

\section{Introduction}

Leavitt path algebras have been introduce by G. Abrams and G. Aranda Pino (see \cite{Abrams}) in 2005, as algebraic analogues of graph C*-algebras. Right after the definition of these algebras there was a spur of activity in the subject, as researches established their structure and found applications to various topics in algebra (see \cite{Abrams},\cite{Abrams1}, \cite{Abrams2}, or \cite{Tomforde}). Two years after the definition of Leavitt path algebras, Mark Tomforde proved the analogue of the graph C*-algebras uniqueness theorems to Leavitt path algebras and established the relation between graph C*-algebras and Leavitt path algebras (see \cite{Tomforde}). We should note that neither the graph C*-algebras nor the Leavitt path algebras results are obviously consequences of the others. Actually it is often the case that analogue results have completely different proofs and, moreover, neither result can be seen to imply the other. Also recently, in \cite{AraBrustenga}, the relations between the theory of quiver representations and the theory of representations of Leavitt path algebras were explored.

It is in the spirit above that we write this paper. Our aim is to prove analogue versions of the representation theorems (for graph C*-algebras) in \cite{repgraph} and \cite{uniteq}, that is, to show how to obtain representations of Leavitt path algebras from E-algebraic branching systems, to study these representations and to give sufficient conditions to guarantee that a representation of $L_K(E)$ is equivalent to a representation arising from an E-algebraic branching system.

As it is often the case, we use many different techniques from the ones used in \cite{repgraph} and \cite{uniteq}, and what is even more interesting, we are able to obtain deeper versions, for Leavitt path algebras, of the results in \cite{repgraph} and \cite{uniteq}. Namely, we are able to state a sufficient condition to guarantee faithfulness of a representation induced by an E-algebraic branching system (We note that our condition is still valid even in the case of a graph $E$ that does not satisfy condition (L), in which case the Cuntz-Krieger Uniqueness theorem of \cite{Tomforde} fail).

Among other things, we expect that the concrete faithful representations of Leavitt path algebras that we present here will deepen, and at the same time make it easier, the understanding of these algebras. Furthermore, we expect our results in the equivalence of representations to be useful for the study of irreducible representations of Leavitt path algebras.

The paper is organized as follows: Below we recall some basic terminology and definitions about Leavitt path algebras, following \cite{Tomforde}. We devote section 2 to the introduction of E-algebraic branching systems and the representations of $L_K(E)$ induced by then. In section 3 we show that, for any graph $E$, we may always find representations induced by E-algebraic branching systems. We present one of the main results of the paper in section 4, where we show that for any graph with no sinks it is possible to construct faithful representations arising from E-algebraic branching systems. In order to do so, we also present a sufficient condition for a representation arising from a E-algebraic branching system to be faithful. In section 5, we make precise what we mean by equivalence of representations, and give a sufficient condition to guarantee that a given representation (or a restriction of it) is equivalent to a representation arising from an E-algebraic branching system. Finally, in section 6, we show that for certain graphs the sufficient condition of section 5 is always satisfied, that is, any representation (or a restriction of it) of $L_K(E)$ is equivalent to a representation arising from an E-algebraic branching system.

Before we proceed, let us recall some definitions:

By a graph we always mean a directed graph $E=(E^0, E^1, r, s)$, where $E^0$ is a countable set of vertices, $E^1$ is a countable set of edges and $r,s:E^1\rightarrow E^0$ are the range and source maps. A path is a sequence $\alpha:= e_1e_2 \ldots e_n$ of edges with $r(e_i)=s(e_{i+1})$, for $1 \leq i \leq n-1$ and we say that the path $\alpha$ has length $|\alpha|:=n$. We denote the set of paths of length $n$ by $E^n$ and consider the vertices in $E^0$ to be paths of length zero. We also let $E^*:= \cup_{n=0}^{\infty} E^n$ denote the paths of finite length.

We let $(E^1)^*$ denote the set of formal symbols $\{e^*:e\in E^1\}$ and for $\alpha:= e_1e_2 \ldots e_n \in E^n$ we define $\alpha^*:= e_n^*e_{n-1}^* \ldots e_1^*$ . We also define $v^*=v$ for all $v\in E^0$.

\begin{definicao}\label{defleviatt}(As in \cite{Tomforde}). Let $E$ be a directed graph, and $K$ be a field. The Leavitt path algebra of $E$ with coefficients in $K$, denoted $L_K(E)$, is the universal $K$-algebra generated by a set $\{v:v\in E^0\}$, of pairwise orthogonal idempotents, together with a set $\{e,e^*:e\in E^1\}$ of elements satisfying:
\end{definicao}
\begin{enumerate}
\item $s(e)e=er(e)=e$ for all $e\in E^1$
\item $r(e)e^*=e^*s(e)=e^*$ for all $e\in E^1$ 
\item $e^*f= \delta_{e,f} r(e)$ for all $e,f \in E^1$
\item $v=\sum\limits_{e\in E^1:s(e)=v}ee^*$ for every vertex $v$ with $0<\#\{e:s(e)=v\}<\infty$

\end{enumerate}

\section{E-algebraic branching systems}

In this section we will define E-algebraic branching systems associated to a directed graph $E$ and we will show how these E-algebraic branching systems induce representations of the associated Leavitt path algebra, in the $K$ algebra of the homomorphisms in a certain module.

We start with the definition of an E-algebraic branching system:

\begin{definicao}\label{brancsystem}
Let $X$ be set and let $\{R_e\}_{e\in E^1}$, $\{D_v\}_{v\in E^0}$ be families of subsets of $X$ such that:
\begin{enumerate}
\item $R_e\cap R_d= \emptyset$ for each $d,e\in E^1$ with $d\neq e$,
\item $D_u\cap D_v= \emptyset$ for each $u,v\in E^0$ with $u\neq v$,
\item $R_e\subseteq D_{s(e)}$ for each $e\in E^1$,
\item $D_v=\bigcup\limits_{e:s(e)=v}R_e$\,\,\,\,\, if\,\,\,\,\, $0<\#\{e\in E^1\,\,:\,\,s(e)=v\}< \infty$,
\item for each $e\in E^1$, there exists a bijective map $f_e:D_{r(e)}\rightarrow R_e$.
\end{enumerate}

A set $X$, with families of subsets $\{R_e\}_{e\in E^1}$, $\{D_v\}_{v\in E^0}$, and maps $f_e$ as above, is called an $E$- algebraic branching system, and we denote it by $(X,\{R_e\}_{e\in E^1}, \{D_v\}_{v\in E^0}, \{f_e\}_{e\in E^1})$, or when no confusion arises, simply by $X$.
\end{definicao}


Next, fix an $E$-algebraic branching system $X$. Let $M$ be the $K$ module of all functions from $X$ taking values in $K$ and let $Hom_K(M)$ denote the $K$ algebra of all homomorphisms from $M$ to $M$ (with multiplication given by composition of homomorphisms and the other operations given in the usual way).

Now, for each $e\in E^1$ and for each $v\in E^0$, we will define homomorphisms $S_e$, $S_e^*$ and $P_v$ in $Hom_K(M)$.

Let $S_e$ be defined as follows: 
$$\left( S_e \phi \right) (x) = 
\begin{cases} \phi(f_e^{-1}(x)),  \text { if }  x\in R_e \\ 0,  \text{ if } x\notin R_e \end{cases}, $$ where $\phi$ is a function in $M$.
 
In order to simplify notation, in what follows we will make a small abuse of the characteristic function symbol and denote the above homomorphism by: 
$$S_e \phi=\chi_{R_e}\cdot \phi\circ f_e^{-1}.$$ 

In a similar fashion to what is done above, and making the same abuse of the characteristic function symbol, we define the homomorphism $S_e^*$ by $$S_e^* \phi=\chi_{D_{r(e)}}\cdot \phi\circ f_e,$$ where $\phi \in M$.

Finally, for each $v\in E^0$, and for $\phi \in M$, we define $P_v $ by $$P_v \phi=\chi_{D_v} \cdot \phi,$$
that is, $P_v$ is the multiplication operator by $\chi_{D_v}$, the characteristic function of $D_v$. 

 \begin{teorema}\label{rep} Let $X$ be an $E$- algebraic branching system. Then there exists a representation (that is, an algebra homomorphism) $\pi: L_K(E)\rightarrow Hom_K(M)$ such that 
$$\pi(e)= S_e, \text{ } \pi(e^*)= S_e^* \text{ and } \pi(v)=P_v,$$ for each $e\in E^1$ and $v\in E^0$.
  \end{teorema}
  
\demo 
Since $L_K(E)$ is an universal object, all we need to do is show that the families $\{S_e, S_e^*\}_{e\in E^1}$ and $\{ P_v \}_{v\in E_0}$ satisfy the relations given in definition \ref{defleviatt}. 


It is clear that all $P_v$ are idempotents, and orthogonality follows from item 2 in definition \ref{brancsystem}.

Now, let $\phi \in M$. Notice that,

$$P_{s(e)} S_e (\phi) = \chi_{D_{s(e)}} \cdot S_e(\phi) =  \chi_{D_{s(e)}} \cdot \chi_{R_e} \cdot \phi \circ f_e^{-1} = S_e(\phi),$$ where the last equality follows from condition 3 in definition \ref{brancsystem}.
In a similar way, one shows that $S_e P_{r(e)} = S_e$ and we have relation 1 in \ref{defleviatt}. Relation 2 of the definition of the Leavitt path algebras follows analogously.

To see that relation 3 holds notice that $$S_e^*S_g ( \phi) =  \chi_{D_{r(e)}} \cdot \left( S_g(\phi) \circ f_e \right) =  \chi_{D_{r(e)}} \cdot \chi_{R_g} \circ f_e \cdot \phi \circ f_g^{-1} \circ f_e = \delta_{e,g} P_{r(e)},$$ where we used that  $R_e\cap R_g= \emptyset$, for $g\neq e$, to obtain the last equality.

Finally, notice that if  $0<\{e\in E^1\,:\, s(e)=v\}<\infty$  then $D_v=\bigcup\limits_{e:s(e)=v}R_e$, and hence $$ \sum\limits_{\{e:s(e)=v\}} S_e S_e^* (\phi) = \sum\limits_{\{e:s(e)=v\}} \chi_{R_e} \cdot S_e^* (\phi) \circ f_e^{-1} = \sum\limits_{\{e:s(e)=v\}} \chi_{R_e} \cdot \chi_{D_{r(e)}}\circ f_e^{-1} \cdot \phi =$$
$$=\sum\limits_{\{e:s(e)=v\}} \chi_{R_e} \cdot \phi=\chi_{D_v}\cdot \phi =P_v(\phi). $$

\fim  

\begin{obs}\label{Mvanishing} Notice that theorem \ref{rep} still holds if we change the module $M$ of all functions from $X$ to $K$ for the module of all functions from $X$ to $K$ that vanish in all, but a finite number of points, of $X$.
\end{obs}

In the next section we consider the question of existence of E-algebraic branching systems(and their induced representations) for any given graph $E$.

\section{Existence of E-algebraic branching systems}

Let $E$ be a graph, with $E^0$ and $E^1$ countable. Next we show that there exists an E-algebraic branching system in $\R$ associated to $E$. Our proof is constructive and one can actually obtain a great number of E-algebraic branching systems following the ideas below.

\begin{teorema}\label{existencebrancsys}
Let $E=(E^0,E^1,r,s)$ be a graph, with $E^0,E^1$ both countable. Then there exists an E-branching  system $X$, where $X$ is an (possible unlimited) interval of $\R$.
\end{teorema}

\demo 
Let $E^1=\{e_i\}_{i=1}^\infty$ (or, if $E^1$ is finite, let $E^1=\{e_i\}_{i=1}^N$). For each $i\geq 1$ define $R_{e_i}=[i-1,i)$. 
Let $W=\{v\in E^0\,:\,\,v \text{ is a sink}\}$ (a vertex $v\in E^0$ is a sink if $v\notin s(E^1)$). Note that $W$ is finite or infinite countable. Write $W=\{v_i\,:\,\,i=1,2,3,...\}$. For each $v_i\in W$, define $D_{v_i}=[-i,-i+1)$. For the vertices $u\in E^0$ which are not sinks, define $D_u=\bigcup\limits_{e_i:s(e_i)=u}R_{e_i}$. Note that items 1-4 from definition \ref{brancsystem} are satisfied. It remains to define functions which satisfy item 5.

Let $\overline{e}\in E^1$. 

If $r(\overline{e})$ is a sink then $r(\overline{e})=v_i\in W$, and so $D_{r(\overline{e})}=[-i, -i+1)$. Then we define $f_{\overline{e}}:D_{r(\overline{e})}\rightarrow R_{\overline{e}}$ as any bijection between these sets (for example, the linear bijection). 

If $r(\overline{e})=\overline{v}$ is not a sink, then $$D_{r(\overline{e})}=D_{\overline{v}}=\bigcup\limits_{e:s(e)=\overline{v}}R_e.$$ 
To define the function $f_{\overline{e}}:D_{r(\overline{e})}\rightarrow R_{\overline{e}}$ in this case we proceed as follows. 
 
 First we divide the interval $R_{\overline{e}}$ in $\#\{e:s(e)=\overline{v}\} $ pairwise disjoint (open on the right and closed on the left hand side) intervals $I_e$ (notice that we might have to divide $R_{\overline{e}}$ in a countable infinite number of intervals). Then, we define $\tilde{f_{\overline{e}}}:\bigcup\limits_{e:s(e)=\overline{v}}R_e\rightarrow \bigcup\limits_{e:s(e)=\overline{v}} I_e$ so that $\tilde{f_{\overline{e}}}_{|_{R_e}}$ is a bijection between $R_e$ and $I_e$ (for example, the linear bijection).

Now, defining $$X=\left(\bigcup\limits_{e_i\in E^1}R_{e_i}\right)\cup\left(\bigcup\limits_{v_i\in W}D_{v_i}\right)$$
 we obtain the desired E-algebraic branching system.
 
\fim

Theorem \ref{existencebrancsys} together with theorem \ref{rep} guarantees that every Leavitt path algebra $L_K(E)$ of a countable graph $E$  may be represented in $Hom_K(M)$. Let us summarize this result in the following corollary:

\begin{corolario}\label{cor1} Given a countable graph $E$, there exists a homomorphism $\pi:L_K(E)\rightarrow Hom_K(M)$ such that $$\pi(v)(\phi)=\chi_{D_v}.\phi,\,\,\,\,\,\pi(e)(\phi)=\chi_{R_e}.\phi\circ f_e^{-1}\,\,\,\,\text{ and }\,\,\,\,\pi(e^*)(\phi)=\chi_{D_{r(e)}}.\phi\circ f_e$$ for each $\phi\in M$,  where $M$ is the $K$ module of all functions from $X$ taking values in $K$, $X$ is an (possible unlimited) interval of $\R$, and $R_e$ and $D_v$ are as in theorem \ref{existencebrancsys}

\end{corolario}

We now seek conditions that guarantee the faithfulness of the representations we have constructed above (of course when the Leavitt path algebra is simple any non-zero representation is faithful).

\section{Faithful representations of Leavitt path algebras of row-finite graphs without sinks}

Important results regarding faithfulness of a representation in the literature include the Graded Uniqueness theorem and the Cuntz-Krieger Uniqueness theorem (see \cite{Tomforde}). In fact, we may use the Cuntz-Krieger Uniqueness theorem for the representations of corollary \ref{cor1}. This theorem guarantees that for any graph $E$ that satisfies condition $(L)$ (each closed path in $E$ has an exit, that is, if $\alpha=\alpha_1...\alpha_n\in E^n$ with $s(\alpha)=r(\alpha)$, then there exists $e\in E$ such that $s(e)=s(\alpha_i)$ and $e\neq e_i$ for some $i$) faithfulness of a representation follows simply by checking that the representation does not vanish at the vertices of the graph. This follows promptly for the representations of corollary \ref{cor1} and hence, for graphs that satisfy condition $(L)$, they are faithful.

As we could see above, the Cuntz-Krieger Uniqueness theorem is a very powerful tool, but it excludes some very simple examples, as for the graph $E$ defined by $E^0=\{ *\}$ and $E^1=\{x\}$ ($E$ consists of one vertex and one "loop" edge). The Leavitt path algebra associated to this graph is $K[x,x^{-1}]$, the Laurent polynomials algebra, see \cite{Abrams}.

In order to overcome problems as the one mentioned above, in this section we introduce a sufficient condition (valid for row finite graphs without sinks) to guarantee that a representation of $L_K(E)$ induced by an E-algebraic branching system is faithful.  

Recall that a graph is row-finite if $s^{-1}(v)$ is finite, for each $v\in E^0$, and a sink is a vertex which emits no edges. 

Let $(X,\{R_e\}_{e\in E^1},\{D_v\}_{v\in E^0}, \{f_e\}_{e\in E^1})$ be an $E$-algebraic branching system. A closed path $\alpha=e_1...e_n$ in the graph $E$ is a path such that $r(e_i)=s(e_{i+1})$ and $r(\alpha):=r(e_n)=v=s(e_1)=:s(\alpha)$. For a closed path $\alpha$, let $$f_\alpha:D_v\rightarrow R_{e_1}\subseteq D_v$$ denote the composition $$f_\alpha:=f_{e_1}\circ...\circ f_{e_n}.$$ 

\begin{obs} Notice that since $\alpha$ is a path $f_\alpha$ is well defined.
\end{obs}

\begin{teorema}\label{faithfulrep} Let $(X,\{R_e\}_{e\in E^1},\{D_v\}_{v\in E^0}, \{f_e\}_{e\in E^1})$ be an $E$-algebraic branching system for a row-finite graph without sinks $E$. Suppose that for each finite set of closed paths $\{\alpha^1,...,\alpha^n\}$ in $E$, beginning on the same vertex $v$, there is an element $z_0\in D_v$ such that $f_{\alpha^i}(z_0)\neq z_0$ for all $i\in \{1,...,n\}$. Then, the representation of $L_K(E)$ induced by this E-algebraic branching system is faithful.
\end{teorema}

\demo Let $\pi:L_K(E)\rightarrow Hom_K(M)$ be the representation induced by the E-algebraic branching system, as in theorem \ref{rep}. (Recall that $M$ is the $K$-module of functions from $X$ to $K$).  

Let $x\in L_K(E)$, $x\neq 0$. Our aim is to show that $\pi(x)\neq 0$. We will separate the proof in a few steps. We start with:

\af{1}{For each $n\in\N$, there exists a path $e_1...e_n$ of length $n$ such that $xe_1...e_n\neq 0$.} 

It is enough to show that for all $0\neq y$ in $L_K(E)$ there is $e\in E^1$ such that $ye\neq 0$.
Since $x$ has a right identity, given by a sum of projections $v\in E^0$, then there is a $v\in E^0$ such that $xv\neq 0$. Since $v=\sum\limits_{e\in s^{-1}(v)}ee^*$ then for some $e$, $xe\neq 0$, otherwise $xv=\sum\limits_{e\in s^{-1}(v)}xee^*=0$.





\fimaf

In order to state our next step we need to make a few observations.

First, notice that, from step 1, we may find a sufficiently large $n$ such that the product $xe_1...e_n$ may be written as a finite sum:

$$xe_1...e_n=\sum\limits_{i=1}^p\gamma_ic^i \neq 0,$$ where $c^i$ are paths in $E$ with $|c^i|\geq 1$ for all $i$, $\gamma_i\neq 0$ for all $i$, and $c^i\neq c^j$ for $i\neq j$.

Also, for each $z\in \bigcup\limits_{u\in E^0} D_u$, denote by $\delta_z$ the function defined by $\delta_z(y)=[y=z]$, where $[y=z]=0$ if $y\neq z$ and $[y=z]=1$ if $y=z$. 

In the next step we characterize how $\pi(e_1...e_n)$ acts on $\delta_z$.

\af{2}{For each path $d_1...d_n$ in $E$ and each $\delta_z$ we have that $$\pi(d_1...d_n)(\delta_z)=[z\in D_{r(d_n)}]\delta_{f_{d_1...d_n}(z)}.$$

The proof of this step follows from the fact that $\pi(e)(\phi)=\chi_{R_e}.\phi\circ f_e^{-1}$ for each edge, and is left to the reader.
\fimaf

Next, fix $c\in \{c^1,...,c^p\}$ such that $|c|\leq |c^i|$ for all $i\in\{1,...,p\}$. Notice that, from Step 2, $\pi(c)(\delta_{z})=\delta_{f_c(z)}$ for all $z\in D_{r(c)}$, and we have:

\af{3}{Let $\alpha$ be a path in $E$. If $r(\alpha)\neq r(c)$ then $\pi(\alpha)(\delta_{z})=0$ for all $z\in D_{r(c)}$}.

To see this, just notice that if $r(c)\neq r(\alpha)$ then $D_{r(\alpha)}\cap D_{r(c)}=\emptyset$. So, for $z\in D_{r(c)}$, by Step 2 it follows that $\pi(\alpha)\delta_z=0$. 
\fimaf

\af{4}{For all $c^i\in\{c^1,...,c^p\}\setminus\{c\}$ with $|c^i|= |c|$ we have that $$\pi(c^i)(\delta_z)(f_c(z))=0$$ for all $z\in D_{r(c)}$.}

The proof of this step goes as follows:

Write $c^i=d_1...d_n$ and $c=c_1...c_n$. Let $j_0$ be the smallest of the indexes $j$ such that $d_j\neq c_j$, and let $z\in D_{r(c)}$. 

We claim that $f_c(z)\neq f_{c^i}(z)$. Suppose not. Then 
$$f_{c_1}\circ....\circ f_{c_{j_0-1}}\circ f_{c_{j_0}}\circ...\circ f_{c_n}(z)=f_{c_1}\circ....\circ f_{c_{j_0-1}}\circ f_{d_{j_0}}\circ...\circ f_{d_n}(z)$$ and since the $f_e$'s are bijections for all $e\in E^1$, the above equality implies that $$f_{c_{j_0}}\circ...\circ f_{c_n}(z)= f_{d_{j_0}}\circ...\circ f_{d_n}(z),$$ which is a contradiction, since $c_{j_0}\neq d_{j_0}$ implies that the image of $f_{c_{j_0}}$ is disjoint from the image of $f_{d_{j_0}}$.

So, $f_c(z)\neq f_{c^i}(z)$, and hence

$$\pi(c^i)(\delta_z)(f_c(z))=\delta_{f_{c^i}(z)}(f_c(z))=0.$$
\fimaf

\af{5}{There exists $z_0\in D_{r(c)}$ such that $\pi(c^i)(\delta_{z_0})(f_c(z_0))$=0 for all $c^i\in \{c^1,...,c^p\}\setminus \{c\}$}.

First notice that, by steps 3 and 4, if $r(c^i)\neq r(c)$, or if $|c^i|=|c|$ and $c^i\neq c$, then $\pi(c^i)(\delta_z)(f_c(z))=0$ for all $z\in D_{r(c)}$.

So, let $c^i$ be such that $|c^i|>|c|$ and $r(c^i)=r(c)$. Write $c=c_1...c_n$ and $c^i=d_1...d_n...d_m$. If $c_j\neq d_j$, for some $j\in \{1,...,n\}$, then proceeding analogously as in the proof of step 4, we have that $f_c(z)\neq f_{c^i}(z)$ for all $z\in D_{r(c)}$, and hence $\pi(c^i)(\delta_z)(f_c(z))=\delta_{f_{c^i}(z)}(f_c(z))=0$ in this case.
 
We are left with the case when $c^i$ is an element of $W$, where $W\subseteq \{c^1,...,c^p\}\setminus \{c\}$ is defined by $W=\{c^i: c^i=c\alpha^i\text{ and } r(c^i)=r(c)\}$. Notice that each $\alpha^i$ is a closed path from $r(c)$ to $r(c)$. By the hypothesis of the theorem, there exists $z_0\in D(r(c))$ such that $f_{\alpha^i}(z_0)\neq z_0$ for all $\alpha^i$. Then, 
$$f_{c^i}(z_0)=f_c\circ f_{\alpha^i}(z_0)\neq f_c(z_0)$$ for all $c^i\in W$ and step 5 is proved.
\fimaf

Let us now conclude the proof of this theorem. Recall that we started with a $x\in L_K(E)$, $x \neq 0$, and have considered the element

$$xe_1...e_n=\sum\limits_{i=1}^p\gamma_ic^i,$$ where $\gamma_i\neq 0$ and $c^i\neq c^j$ for $i\neq j$. 

By the previous steps, there exists a $c\in \{c^1,...,c^p\}$ and a $z_0\in D_{r(c)}$ such that $\pi(c^i)(\delta_{z_0})(f_c(z_0))=0$ for all $c^i\in \{c^1,...,c^p\}\setminus\{c\}$.
But then, 

$$(\pi(x)\pi(e_1...e_n))(\delta_{z_0})(f_c(z_0))=\pi(xe_1...e_n)(\delta_{z_0})(f_c(z_0))=$$
$$=\sum\limits_{i=1}^p\gamma_i\pi(c^i)(\delta_{z_0})(f_c(z_0))=\gamma_{i_0}\delta_{f_c(z_0)}(f_c(z_0))=\gamma_{i_0},$$

where $i_0$ is such that $c=c^{i_0}$. 

So, it follows that $\pi(x)\pi(e_1...e_n)\neq 0$, and hence $\pi(x)\neq 0$.
\fim

Theorem \ref{faithfulrep} is an important result. It allow us to construct faithful representations of $L_K(E)$, when $E$ is a row finite graph without sinks. For these graphs we describe E-algebraic branching systems that satisfy the conditions of theorem \ref{faithfulrep} (and hence induce faithful representations) below.

So, let $E$ be a countable row-finite graph without sinks.

Since $E^1$ is finite or infinite countable we have that $E^1=\{e_1,e_2,e_3,...,e_N\}$ or $E^1=\{e_1,e_2,e_3,...\}$. 
For each $e_i\in E^1$ define $R_{e_i}:=[i-1,i)$. For a vertex $v\in E^0$ define $D_v:=\bigcup\limits_{e\in s^{-1}(v)}R_e$. 
We also need to define bijective maps $f_e:D_{r(e)}\rightarrow R_e$, for each $e\in E^1$. To do this, first fix an irrational number $\theta\in [0,1)$, and let $h_\theta:[0,1)\rightarrow [0,1)$ be defined by $h_\theta(x)=(x+\theta)mod(1)$, which is a bijective map. For $a,b\in \R$ with $a<b$, define $g_a^b:[0,1)\rightarrow [a,b)$ by $g_a^b(x)=bx+(1-x)a$, which is also a bijective map, with inverse $(g_a^b)^{-1}$.

Let $e\in E^1$. 





Since $r(e)$ is not a sink and $E$ is row-finite then $0<\#\{s^{-1}(r(e))\}<\infty$. So, $s^{-1}(r(e))=\{e_{i_1},...,e_{i_P}\}$ and hence

$$D_{r(e)}=\bigcup\limits_{k=1}^P R_{e_{i_k}}.$$ Since $e=e_j$, for some $j\in \N$, we have that $R_{e}=[j-1,j)$. Write  $R_e$ as the following disjoint union:

$$R_{e}=[j-1,j)=\bigcup\limits_{k=1}^{P}\left[j-1+\frac{k-1}{P},j-1+\frac{k}{P}\right).$$

Now, given $x\in D_{r(e)}$, we have that $x\in R_{e_{i_k}}=[i_k-1,i_k)$ for some $k\in \{1,...,P\}$, and we define 

$$f_{e}(x):=\left(g_{(j-1+\frac{k-1}{P})}^{(j-1+\frac{k}{P})}\circ h_\theta\circ\left(g_{(i_k-1)}^{i_k}\right)^{-1}\right)(x), $$
that is, $f_{e}$ restricted to $R_{e_{i_k}}$ is the composition

$$R_{e_{i_k}}\stackrel{\left(g_{(i_k-1)}^{i_k}\right)^{-1}}{\longrightarrow}[0,1)\stackrel{h_\theta}{\longrightarrow}[0,1)\stackrel{g_{(j-1+\frac{k-1}{P})}^{(j-1+\frac{k}{P})}}{\longrightarrow}\left[j-1+\frac{k-1}{P},j-1+\frac{k}{P}\right).$$

This defines $f_{e}:D_{r(e)}\rightarrow R_{e}$ as a bijective map, and it is not hard to see that $f_{e}(x)=\frac{x+\theta+r_{e}(x)}{P}$, where $r_{e}(x)$ is a rational number, for each $x\in D_r(e)$.

So, for each $e\in E^1$, we have defined a bijective map $f_e:D_{r(e)}\rightarrow R_e$, such that $$f_e(x)=\frac{x+\theta+r_e(x)}{P_e},$$ where $r_e(x)$ is a rational number, for each $x\in D_{r(e)}$ and $P_e$ is a natural number, namely $P_e=\#\{s^{-1}(r(e))\}$.

Defining $X=\bigcup\limits_{e\in E^1}R_e$ we obtain an $E$-algebraic branching system $$(X, \{D_u\}_{u\in E^0}, \{R_e\}_{e\in E^1}, \{f_e\}_{e\in E^1}),$$ and hence we obtain a representation $\pi:L_K(E)\rightarrow Hom_K(M)$ (as in theorem \ref{rep}). 

\begin{corolario}\label{faith} Let $E$ be a row finite graph with no sinks. Then the representation $\pi:L_K(E)\rightarrow Hom_K(M)$ induced by the $E$-algebraic branching system constructed above is faithful.
\end{corolario}

\demo All we need to do is verify the hypothesis of theorem \ref{faithfulrep}, that is, we need to check that for each finite set $\{\alpha^1,...,\alpha^N\}$ of closed paths beginning on the same vertex $v$, there exists an element $z_0\in D_v$ such that $f_{\alpha^i}(z_0)\neq z_0$ for all $i\in \{1,...,N\}$. 

So, let $\alpha=c_1...c_n$ be a closed path beginning on $v$. 
Notice that, for each $x\in D_v$
$$f_{c_1}\circ...\circ f_{c_n}(x)=\frac{x}{P_{c_1}...P_{c_n}}+\theta\left(\frac{1}{P_{c_1}...P_{c_n}}+\frac{1}{P_{c_1}...P_{c_{n-1}}}+...+\frac{1}{P_{c_1}}\right)+\overline{r}(x),$$ where $\overline{r}(x)$ is a rational number and $P_{c_1},..,P_{c_n}$ are natural numbers. It follows that, if $x\in D_v$ is a rational number, then $f_{c_1}\circ...\circ f_{c_n}(x)$ is a irrational number and hence no rational number is a fixed point for $f_{\alpha}$. 
Then, for any finite set $\{\alpha^1,...,\alpha^N\}$ of closed paths in $E$ beginning on $v$, we may choose $z_0\in D_v$ to be a rational number, and so $f_{\alpha^i}(z_0)\neq z_0$ for all $i\in \{1,...,N\}$ as desired.

\fim

\begin{exemplo} Let $E^0 = \{*\}$, $E^1=\{x\}$ as in the figure below.
\end{exemplo}

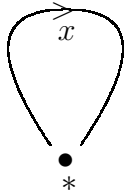
\centerline{
\setlength{\unitlength}{2cm}
\begin{picture}(0,1)
\put(0,0){\circle*{0.08}}
\put(-0.03,-0.2){$*$}
\qbezier(-0.1,0.1)(-0.7,1)(0,1)
\qbezier(0.1,0.1)(0.7,1)(0,1)
\put(-0.1,0.94){$>$}
\put(-0.05,0.8){$x$}
\end{picture}}
\vspace{0.6cm}
 
Notice that $L_K(E) = K[x,x^{-1}]$, the Laurent polynomials in $x$ and $x^{-1}$. By corollary \ref{faith} above, the representation induced by the E-algebraic branching system $X$, where $R_e=[0,1]$, $D_{r(e)}=[0,1]$ and $f_e:D_{r(e)} \rightarrow R_e$ is defined by $f_e(x) = x + \theta \text{ mod } 1$ (that is, $f_e$ is rotation by an irrational number $\theta$) is faithful. It follows that the $K$ algebra of the Laurent polynomials in $x$ and $x^{-1}$ is isomorphic to the sub algebra of $Hom_K(M)$ generated by $\{ S_e, S_e^*\}$, where, for $f\in M$, $S_e \phi = \phi \circ f_e^{-1}$ and $S_e^{*} \phi = \phi \circ f_e$.

\section{Equivalence of representations of $L_K(E)$}

In the previous sections, we have introduced a class of representations of the Leavitt path algebras induced by E-algebraic branching systems. One question which remains is if any representation may be obtained in such a manner. 

In this section, we show that under a certain condition over a graph $E$, each $K$-algebra homomorphism $\widetilde{\pi}: L_K(E)\rightarrow A$ has a sub-representation associated to it which is equivalent to a representation induced by an $E$-algebraic branching system. This is not true in general and we will make a more precise argument just after definition \ref{equivrepres}, where equivalence of representations if formally defined.

Before we proceed, notice that given a $K$-algebra $A$, there exist a $K$-module $V$ and an injective $K$-algebra homomorphism $\varphi:A\rightarrow Hom_K(V)$.  
To see this, note that $A\times K$ is an unital $K$-algebra, with the operations defined by $(a,k)+(b,l):=(a+b,k+l)$, $k(a,l):=(ka,kl)$ and $(a,k)(b,l):=(ab+la+kb,kl)$ for each $a,b\in A$ and $k,l\in K$.
 In particular, $V:=A\times K$ is a $K$-module. Defining, for each $a\in A$, $\varphi(a):V\rightarrow V$ by $\varphi(a)(b,k)=(a,0)(b,k)$, we obtain a injective homomorphism $\varphi:A\rightarrow Hom_K(V)$.
 
Given a $K$-algebra homomorphism $\widetilde{\pi}:L_K(E)\rightarrow A$, using the previous injective $K$-algebra homomorphism $\varphi$, we may consider the composition $\varphi\circ\widetilde{\pi}:L_K(E)\rightarrow Hom_K(V)$. With this in mind, from now on, we only consider representations ($K$-algebra homomorphism) from $L_K(E)$ to $Hom_K(V)$, where $V$ is a $K$-module.  

Next we will prepare the ground for the results in this section. We start with a representation $\Phi:L_K(E)\rightarrow Hom_K(V)$ and define $K$-submodules$$V_u=\Phi(u)(V)$$ and $$V_e=\Phi(e)\Phi(e^*)(V),$$ for all $u\in E^0$ and all $e\in E^1$.
 Since $\Phi$ is a representation of $L_K(E)$, it satisfies the relations of Definition $\ref{defleviatt}$, and it follows that:
 \begin{enumerate}
 \item $V_e\subseteq V_{s(e)}$ for each $e\in E^1$,
 \item $V_e\cap V_f=0$ for each $e,f\in E^1$, $e\neq f$,
 \item $V_u\cap V_w=0$ for each $u,w\in E^0$, $u\neq w$,
 \item $\Phi(e):V_{r(e)}\rightarrow V_e$ is a $K$-module isomorphism, with inverse $\Phi(e^*)$,
 \item $V_u=\bigoplus\limits_{e:s(e)=u}V_e$ if $0<\#\{e:s(e)=u\}<\infty$
 \item $V_u=\left(\bigoplus\limits_{e:s(e)=u}V_e\right)\bigoplus\overline{V_u}$ if $\#\{e:s(e)=u\}=\infty$, where $\overline{V_u}$ is some $K$-submodule of $V_u$,
 \item $V=\left(\bigoplus\limits_{u\in E^0}V_u\right)\bigoplus \overline{V}$, where $\overline{V}$ is a $K$-submodule of $V$. 
\end{enumerate}

To obtain the equality of item 6 above, notice that $V_u$ is a $K$-vector space, and hence we may complete the (Hammel) basis of $\bigoplus\limits_{e:s(e)=u}V_e$ to obtain a basis of $V_u$. The same holds for the last equality.
  
We now intend to pick a particular basis for the K-vector space $V$. By equality 7 above, we need to choose a basis for $V_u$, $u\in E^0$, and $\overline{V}$. 

Before picking the basis for $V_u$, notice that, since $V_e$ and $\overline{V_u}$ are $K$-vector spaces, there exists Hammel basis $\{m_x:x\in R_e\}$ for each $V_e$ and $\{m_x:x\in \overline{I_u}\}$ for each $\overline{V_u}$. 
Choose the index sets $R_e$ and $\overline{I_u}$ as being pairwise disjoint, that is, $R_e\cap R_f=\emptyset$, $R_e\cap \overline{I_u}=\emptyset$ and $\overline{I_u}\cap\overline{I_w}=\emptyset$. Now, we define the basis of $V_u$ in the following way:

\begin{itemize}
\item if $u\notin s^{-1}(E^1)$, choose some basis $\{m_x:x\in D_u\}$ of $V_u$, where $D_u$ is an index set of the basis.
\item if $0<\#\{e\in E^1:s(e)=u\}<\infty$ let $D_u:=\bigcup\limits_{R_e:s(e)=u}R_e$ and so $\{m_x:x\in D_u\}$ is a basis of $V_u$.
\item if $\#\{e\in E^1:s(e)=u\}=\infty$ let $D_u:=\left(\bigcup\limits_{R_e:s(e)=u}R_e\right)\cup \overline{I_u}$ and so $\{m_x:x\in D_u\}$ is a basis of $V_u$.
\end{itemize}

The index sets $D_u$ obtained in the second and last items above are pairwise disjoint. 
In the first item, choose the index sets $D_u$ such that the sets $\{D_u\}_{u\in E^0}$ are pairwise disjoint.

Finally, choose a basis $\{m_x:x\in \overline{I}\}$ of $\overline{V}$ and an index set $\overline{I}$ such that $\overline{I}\cap D_u=\emptyset$ for all $u\in E^0$. We have now chosen a basis for $V$.

Let $$W=\bigoplus\limits_{u\in E^0}V_u.$$ Recall that $V=W\bigoplus \overline{V}$. Let $P_1:V\rightarrow W$ and $P_2:V\rightarrow \overline{V}$ be the two canonical projections and $i_1:W \rightarrow V$ and $i_2:\overline{V}\rightarrow V$ be the two canonical inclusions. Notice that, for each $a\in L_K(E)$, it holds that 

$$\Phi(a)=(P_1+P_2)\Phi(a)(i_1\oplus i_2)=P_1\Phi(a)i_1+P_1\Phi(a)i_2=P_1 \Phi(a) ,$$ since $P_2\Phi(a)=0$. So, we will consider the "restriction of $\Phi$ to $W$", that is, the map 

$$\funcao{\Phi_1}{L_K(E)}{Hom_K(W)}{a}{P_1\Phi(a)i_1=\Phi(a)i_1},$$ which
is a representation.

Our aim is to show (under some additional hypothesis) that the representation $\Phi_1$ is equivalent, in some sense, to a representation induced by an $E$-algebraic branching system. So, we need to define the desired branching system.

Let $X=\bigcup\limits_{u\in E^0}D_u$. By the definition of $D_v$ and $R_e$, it is clear that conditions 1-4 of definition \ref{defleviatt} are satisfied. To obtain an $E$-algebraic branching system, we need to define bijective maps $f_e:D_{r(e)}\rightarrow R_e$. Recall that the restriction $\Phi(e):V_{r(e)}\rightarrow V_e$ is a $K$-module isomorphism, with inverse $\Phi(e^*)$, and the sets $D_{r(e)}$ and $R_e$ are the index sets of the basis of $V_{r(e)}$ and $V_e$, respectively. So, if the basis of $D_{r(e)}$ is taken to the basis of $R_e$, that is, if for each $x\in D_{r(e)}$ we have that $\Phi(e)(m_x)=m_y$ for some $y\in R_e$, then the map $D_{r(e)}\ni x\mapsto y\in R_e$ defines a bijective map $f_e$.

So, from now on we assume this additional hypothesis, that is, we assume that: 
\begin{equation}\label{h1}\tag{$B2B$}
\Phi(e)(\{m_x:x\in D_{r(e)}\})=\{m_y:y\in R_e\} ,\text{ for each } e\in E^1,
\end{equation}
which we call condition (\ref{h1}). Notice that condition (\ref{h1}) is equivalent to say that $\Phi(e^*)(\{m_y:y\in R_e\})=\{m_x:x\in D_{r(e)}\}$ for each $e\in E^1$.

We may now define $f_e:D_{r(e)}\rightarrow R_e$ by $f_e(x)=y$, where $y$ is such that $\Phi(m_x)=m_y$.

Notice that the map $f_e$ is bijective, for each $e\in E^1$, and hence the set $X$ with the families $\{R_e\}_{e\in E^1}$, $\{D_u\}_{u\in E^0}$ and $\{f_e\}_{e\in E^1}$ is an $E$-algebraic branching system.

Before we state our next theorem, we need the following definition:

\begin{definicao}\label{equivrepres} Let $\pi: L_K(E) \rightarrow Hom_K(M)$ and $\Phi: L_K(E) \rightarrow Hom_K(W)$ be representations of $L_K(E)$, where $M$ and $W$ are $K$-modules. We say that $\pi$ is equivalent to $\Phi$ if there exists a $K$-module isomorphism $U:W\rightarrow M$ such that the diagram 
\begin{displaymath}
\xymatrix{
W \ar[r]^{\Phi(a)} \ar[d]_{U} &
W \ar[d]^{U} \\
M \ar[r]_{\pi(a)} & M }
\end{displaymath}

commutes, for each $a\in L_K(E)$.

\end{definicao}

\begin{obs} It is not true in general that every representation of $L_K(E)$ is equivalent to a representation induced by an $E$-algebraic branching system. For instance, let $L_K(E)=K[x,x^{-1}]$, the Laurent polynomials. Let $\Phi:L_K(E)\rightarrow Hom_K(K^n)$ be a representation such that $\Phi(e)$ is a $K$-isomorphism which is not a permutation (that is, $\Phi(e)$ is an invertible matrix in $M_n(K)$ which is not a permutation matrix). Notice that if $\pi:L_K(E)\rightarrow Hom_K(M)$ is a representation induced by an $E$-algebraic branching system (with $M$ and $K^n$ being $K$-isomorphic) then $\pi(e)$ is, loosely speaking, a permutation matrix. So, $\pi$ and $\Phi$ are not equivalent representations.    
\end{obs}

For what follows, let $M:=\{g:X\rightarrow K: g(x)\neq 0 \text{ only for finitely many } x\in X\}$, $Y=X\cup \overline{I}$ (recall that $\overline{I}$ is the index set of $\overline{V}$) and $N=\{g:Y\rightarrow K:g(x)\neq 0\text{ only for finitely many }x\in Y\}$. Recall that $W=\bigoplus\limits_{u\in E^0}V_u$. We are now ready to prove the next theorem.

\begin{teorema}\label{equivrep} Let $\Phi:L_K(E)\rightarrow Hom_K(V)$ be a representation. Choose a basis of (the $K$-vector space) $V$ as constructed above. Suppose that this basis satisfies condition (\ref{h1}). Suppose also that $\Phi(e^*)(\overline{V_{s(e)}})=0$, for all $e\in E^1$, where  $\overline{V_{s(e)}}$ was defined in item 6 above. Then: 
\begin{enumerate}
\item There exists a representation $\pi:L_K(E)\rightarrow Hom_K(M)$, induced by an $E$-algebraic branching system, which is equivalent to $\Phi_1$ (the restriction of $\Phi$ to $W$).

\item if $\overline{V}$ (as in item 7 above) may be chosen such that $\Phi(u)(\overline{V})=0$, for each $u\in E^0$, then there exists an $E$-algebraic branching system which induces a representation $\pi:L_K(E)\rightarrow Hom_K(N)$ which is equivalent to $\Phi$.\end{enumerate} 
\end{teorema}

\demo We begin by proving the first part. Let $\left( X, \{R_e\}_{e\in E^1},\{D_u\}_{u\in E^0}, \{f_e\}_{e\in E^1}\right)$ and $M$ be as defined in the paragraphs preceding this theorem. By theorem \ref{rep}, there exists a representation $\pi:L_K(E)\rightarrow Hom_K(M)$ such that $\pi(e)(g)=\chi_{R_e}.g\circ f_e^{-1}$, $\pi(e^*)(g)=\chi_{D_{r(e)}}.g\circ f_e$ and $\pi(u)(g)=\chi_{D_u}.g$. 

Notice that $M$ is a $K$-module with basis $\{\delta_x\}_{x\in X}$, where $\delta_x:X\rightarrow K$ is defined by $\delta_x(y)=0$, if $y\neq x$ and $\delta_x(y)=1$, if $y=x.$

Recall that $\{m_x:x\in X\}$ is a basis of $V$. So, the map $\{m_x:x\in X\}\ni m_x \mapsto \delta_x\in M$ induces a $K$-module isomorphism $U:W\rightarrow M$.

Next we show that $\Phi_1(a)=U^{-1}\circ\pi(a)\circ U$ for each $a\in L_K(E)$. Notice that it is enough to show that $\Phi_1(e)=U^{-1}\circ\pi(e)\circ U$ and $\Phi_1(e^*)=U^{-1}\circ\pi(e^*)\circ U$, for each $e\in E^1$ and $\Phi_1(v)=U^{-1}\circ\pi(v)\circ U$, for all $v\in E^0$. We will verify the second equality and the other two are left to the reader. 

Notice that, for $x\in X$, $\pi(e^*)(\delta_x)=[x\in R_e]\delta_{f_e^{-1}}(x)$ (where $[x\in R_e]=1,$ if $x\in R_e$ and $[x\in R_e]=0,$ if $x\notin R_e$). So, it follows that $U^{-1}\circ\pi(e^*)\circ U(m_x)=[x\in R_e]\delta_x$. 

We now evaluate $\Phi_1(e^*)(m_x)$. If $m_x\in V_u$, with $u\neq s(e)$, then $\Phi_1(e^*)(m_x)=\Phi(e^*)\Phi(s(e))\Phi(u)(m_x)=0$, since $\Phi(s(e))\Phi(u)=0$. Let $m_x\in V_{s(e)}$. Recall that
$$ V_{s(e)}=\left(\bigoplus\limits_{d\in E^1:s(d)=s(e)}V_d\right)\bigoplus\overline{V_{s(e)}}.$$
If $m_x\in V_d$, for some $d\neq e$, then $\Phi_1(e^*)(m_x)=\Phi(e^*)\Phi(d)\Phi(d^*)(m_x)=0$ (since $\Phi(e^*)\Phi(d))=0)$. If $m_x\in \overline{V_{s(e)}}$, then $\Phi(e^*)(m_x)=0$ by hypothesis. It remains to evaluate $\Phi_1(e^*)(m_x)$ for $m_x\in V_e$. In this case, $\Phi_1(e^*)(m_x)=\Phi(e^*)(m_x)=m_{f_e^{-1}}(x)$, by the definition of the map $f_e^{-1}$.  

So, it follows that $U^{-1}\circ\pi(e^*)\circ U=\Phi_1(e^*)$ as desired and we have that
$$\Phi_1(a)=U^{-1}\circ\pi(a)\circ U,$$ for all $a\in L_K(E)$.

Defining $T:Hom_K(W)\rightarrow Hom_K(M)$ by $T(A)=U\circ A\circ U^{-1}$, which is a $K$-algebra isomorphism, we obtain that $\Phi_1$ is equivalent to $\pi$.

To prove the second part of the theorem consider the $E$-algebraic branching system $Y=X\cup \overline{I}$ (recall that $\overline{I}$ is the index set of $\overline{V}$). Consider $\{R_e\}_{e\in E^1}$, $\{D_u\}_{u\in E^0}$ and $\{f_e\}_{e\in E ^1}$  as in the first part. This $E$-algebraic branching system induces a representation $\pi:L_K(E)\rightarrow Hom_K(N)$, where 
$$N=\{g:Y\rightarrow K:g(x)\neq 0\text{ only for finitely many }x\in Y\}.$$
The map $V\ni m_x\rightarrow \delta_x\in N$ induces a $K$-module isomorphism $Q:V\rightarrow N$, and the map $L:Hom_K(V)\rightarrow Hom_K(N)$ defined by $L(A)=Q\circ A\circ Q^{-1}$, for each $A\in Hom_K(V)$ is also an isomorphism. The rest of the proof follows analogously to what was done above for the first part of the theorem.
\fim   

\begin{obs} a) If the graph $E$ is row-finite, then $V_{s(e)}=\bigoplus\limits_{d\in E^1:s(d)=s(e)}V_d$, and the condition $\Phi(e^*)(\overline{V_{s(e)}})=0$ (which appears in the hypothesis of the previous theorem) is vacuously satisfied. So, the first part of the previous theorem applies to any representation of row-finite graphs, as long as (\ref{h1}) is satisfied.

b) If $E^0$ is finite then $\overline{V}$ may be chosen so that $\Phi(u)(\overline{V})=0$, for each $u\in E^0$. In fact, if $E^0$ is finite, define 
$$\overline{V}=\{m\in V: \Phi(u)(m)=0\,\,\,\forall\,\, u\in E^0\}.$$
Then it is clear that $\left(\bigoplus\limits_{u\in E^0}V_u\right)\bigoplus \overline{V}\subseteq V$. For a given $m\in V$, write $m$ as  the sum $$m=\left(\sum\limits_{u\in E^0}\Phi(u)(m)\right)+\left(m-\sum\limits_{u\in E^0}\Phi(u)(m)\right),$$ and note that $\sum\limits_{u\in E^0}\Phi(u)(m)\in\bigoplus\limits_{u\in E^0}V_u$ and $\left(m-\sum\limits_{u\in E^0}\Phi(u)(m)\right)\in \overline{V}$. So, it follows that $$\left(\bigoplus\limits_{u\in E^0}V_u\right)\bigoplus \overline{V}= V,$$ and the second part of the previous theorem applies to representations of any graph $E$, with $E^0$ finite (as long as (\ref{h1}) is satisfied).
\end{obs}

Given a representation $\Phi:L_K(E)\rightarrow Hom_K(V)$, to see if this representation (or its restriction to $W$, $\Phi_1$) is equivalent to a representation induced by an $E$-algebraic branching system, we must be able to, among other things, guarantee the existence of a basis of $V$ satisfying the hypothesis of theorem \ref{equivrep}. The existence of such a basis seems to be intrinsic to the representation $\Phi$ and to the module $V$, however, we prove in the next section that under a certain (sufficient but not necessary) condition over the graph $E$, it is always possible to choose such a basis of $V$.

\section{A sufficient condition over $E$ to guarantee equivalence of representations}

Most of this section is inspired by corresponding results and ideas for graph C*-algebras, as done in \cite{repgraph} and \cite{uniteq}. For the reader's convenience, we adapt the necessary definitions and results below.

\begin{definicao}\cite{uniteq} Let $E=(E^0,E^1,r,s)$ be a graph. We say that:
\begin{enumerate}
\item A path without orientation between $u,v\in E^0$ is a pair of sequences ($u_0u_1...u_n;e_1...e_n$) of vertices $u_i$ and edges $e_j$ such that: $u=u_0$, $v=u_n$, $e_i\neq e_j$ for $i\neq j$, and for each $i$ it holds that $s(e_i)=u_{i-1}$ and $r(e_i)=u_i$, or $r(e_i)=u_{i-1}$ and $s(e_i)=u_i$.
\item A graph $E$ is $P$-simple if for each $u,v \in r(E^1)\cup s(E^1)$, with $u\neq v$, there exists at most one path without orientation between $u$ and $v$, and moreover it does not exist $e\in  E^1$ such that $r(e)=s(e)$.
\item Let $E=(E^0,E^1, r,s)$ be a graph. We say that a subset $Z$ of $E^0$ is connected if, for each $u,v\in Z$, there exists a path without orientation between $u$ and $v$. 

\end{enumerate}
\end{definicao}


For a given graph $E$, $E^0$ is obviously not necessarily connected, but it is always possible to write 

$$E^0=\left(\bigcup\limits_{i\in \Delta}^.Z_i\right)\bigcup\limits^.R,$$
where each $Z_i$ is connected and $R$ is the set of isolated vertices.

\begin{definicao}\cite{uniteq} A vertex $v\in E^0$ is an extreme vertex of $E$ if $\#\{r^{-1}(v)\cup s^{-1}(v)\}=1$ and if there does not exist an edge $e\in E^1$ such that $r(e)=v=s(e)$. If $v$ is an extreme vertex, then the unique edge adjacent to $v$ is called an extreme edge.
\end{definicao} 

We denote by $X_1$ the set of extreme vertices of $E$ (the level 1 vertices) and by $Y_1$ the set of extreme edges of $E$ (the level 1 edges). Notice that $\E_1=(E^1\setminus Y_1,E^0\setminus X_1,r,s)$  is a new graph (here $r$ and $s$ are the restriction maps $r,s:E^1\setminus Y_1\rightarrow E^0\setminus X_1$). We denote by $X_2$ the set of extreme vertices of $\E_1$ (the level 2 vertices), and by $Y_2$ the extreme edges of $\E_1$ (the level 2 edges). Proceeding inductively we define the level $n$ vertices set, $X_n$, and the level $n$ edges set, $Y_n$, if such vertices and edges exist. For more details see $[\ref{uniteq}]$.

Our aim in this section is to describe a sufficient condition, over the graph $E$, which guarantees that a representation of $L_K(E)$ is equivalent to a representation induced by an $E$-algebraic branching system. So, let us fix a representation $\Phi:L_K(E)\rightarrow Hom_K(V)$. Let $V_e:=\Phi(e)\Phi(e^*)(V)$ and $V_u:=\Phi(u)(V)$ for each $e\in E^1$ and $u\in E^0$. By theorem \ref{equivrep}, all we need to do is verify the existence of basis of $V_e$ and $V_u$ satisfying the conditions of that theorem. Next we show the existence of such basis, under a certain condition over the graph $E$.

\begin{teorema}\label{theorem1} Let $E=(E^0,E^1,r,s)$ be a graph such that $r(E^1)\cup s(E^1)$ is connected and suppose 
$$r(E^1)\cup s(E^1)=\bigcup\limits_{n=1}^m X_n \,\,\,\,\,\,\,\text{ or }\,\,\,\,\,\,\, r(E^1)\cup s(E^1)=\bigcup\limits_{n=1}^m X_n\cup\{\overline{v}\}.$$ Let $\Phi:L_K(E)\rightarrow Hom_K(V)$ be a representation. For each $e\in E^1$ and $v\in E^0$, consider the subspaces $V_e:=\Phi(e)\Phi(e^*)(V)$ and $V_u:=\Phi(u)(V)$. Then, there exists basis $B_e$ of $V_e$ and $B_u$ of $V_u$ such that:

\begin{enumerate}
\item[1)] if $e\in s^{-1}(u)$, then $B_e\subseteq B_u$ and if $0<|s^{-1}(u)|<\infty$, then $B_u=\bigcup\limits_{e\in s^{-1}(u)}B_e$;
\item[2)] if $e\in r^{-1}(u)$, then $\Phi(e)(B_u)=B_e$. (and hence the basis satisfies hypothesis (\ref{h1})).
\end{enumerate}
\end{teorema}

\demo 
The proof is similar to the proof of [\ref{uniteq}:Theorem 4.1].
\fim

In the following proposition, we obtain a sufficient condition over the graph $E$ to conclude that  

$$r(E^1)\cup s(E^1)=\bigcup\limits_{n=1}^m X_n \,\,\,\,\,\,\,\text{ or }\,\,\,\,\,\,\, r(E^1)\cup s(E^1)=\bigcup\limits_{n=1}^m X_n\cup\{\overline{v}\}.$$

\begin{proposicao}\label{prop1}\cite{uniteq}
If $r(E^1)\cup s(E^1)$ is finite and connected and if $E$ is P-simple then $r(E^1)\cup s(E^1)=\bigcup\limits_{n=1}^mX_n$ or $r(E^1)\cup s(E^1)=\left(\bigcup\limits_{n=1}^mX_n\right)\bigcup\limits^\cdot \{\overline{v}\}$.
\end{proposicao}

We are now ready to prove the main result of this section, and finally present the reader with the condition over the graph $E$ that guarantees that a representation of $L_K(E)$ is equivalent to a representation induced by an $E$-algebraic branching system.

\begin{teorema}\label{teorONequivrepres}
Let $E=(E^0,E^1,r,s)$ be a graph. Write $$E^0=\left(\bigcup\limits_{i\in \Delta}^.Z_i\right)\bigcup\limits^.R,$$
where each $Z_i$ is connected and $R$ is the set of isolated vertices. 

Suppose $Z_i=\bigcup\limits_{n=1}^{m_i}X_n$ or $Z_i=\bigcup\limits_{n=1}^{m_i}X_n\cup \{\overline{v_i}\}$ for each $i\in \Delta$ (for example, if each graph $E_i:=(r^{-1}(Z_i)\cup s^{-1}(Z_i),Z_i,r,s)$ is $P$-simple and $Z_i$ is finite, see proposition \ref{prop1}).

Let $\Phi:L_K(E)\rightarrow Hom_K(V)$ be a representation and suppose $\Phi(e^*)(\overline{V_e})=0$ for each $e\in E^1$. Then:
\begin{enumerate} 
\item The representation $\Phi_1$ (the restriction of $\Phi$ to $W$) is equivalent to a representation induced by an $E$-algebraic branching system.
\item If $\Phi(u)(\overline{V})=0$, for each $u\in E ^0$, then $\Phi$ is equivalent to a representation induced by an $E$-algebraic branching system.
\end{enumerate} 
\end{teorema}

\demo Let $\Phi:L_K(E)\rightarrow Hom_K(V)$ be a representation. Define $V_e=\pi(e)\pi(e^*)(V)$ and $V_u=\pi(u)(V)$. Note that $V=\bigoplus\limits_{i\in \Delta}\left(\bigoplus\limits_{u\in Z_i}V_u\right)\bigoplus\overline{V}$, where $\overline{V}$ is some submodule of $V$. Applying theorems \ref{theorem1} and theorem \ref{equivrep} to each graph $E_i:=(r^{-1}(Z_i)\cup s^{-1}(Z_i),Z_i,r,s)$, we obtain E-algebraic branching systems:
$$\left(X_i, \{D_v\}_{v\in Z_i}, \{R_e\}_{e\in s^{-1}(Z_i)\cup r^{-1}(Z_i)}, \{f_e\}_{e\in s^{-1}(Z_i)\cup r^{-1}(Z_i)}\right)$$ and $K$-module isomorphisms $U_i:\bigoplus\limits_{u\in Z_i}V_u\rightarrow M_i$, for each $i\in \Delta$. 

Now, the first part of the theorem follows if we consider the representation induced by the $E$-algebraic branching system 
$$\left(\bigcup\limits_{i\in \Delta}X_i, \{R_e\}_{e\in E^1}, \{D_u\}_{u\in E^0},\{f_e\}_{e\in E^0}\right),$$ and the $K$-module isomorphism
$$U:\bigoplus\limits_{i\in \Delta}\left(\bigoplus\limits_{u\in Z_i}V_u\right)\rightarrow \bigoplus\limits_{i\in \Delta}M_i \text{ given by } U:=\bigoplus\limits_{i\in \Delta}U_i.$$

The second statement of the theorem follows if we consider the $E$-algebraic branching system 
$$\left(\left(\bigcup\limits_{i\in \Delta}X_i\right)\cup I, \{R_e\}_{e\in E^1}, \{D_u\}_{u\in E^0},\{f_e\}_{e\in E^0}\right),$$ where $I$ is some index set (with $I\cap X_i=\emptyset$ for each $i\in \Delta$) of a (Hammel) basis $\{m_x:x\in I\}$ of $\overline{V}$, and the $K$-module isomorphism
$$Q:V=\bigoplus\limits_{i\in \Delta}\left(\bigoplus_{u\in Z_i}V_u\right)\bigoplus\overline{V}\rightarrow N$$ defined by $Q(m)=U(m)$, if $m\in \bigoplus\limits_{i\in \Delta}\left(\bigoplus\limits_{u\in Z_i}V_u\right)$, and $Q(m_x)=\delta_x$, if $m_x$ is a element of the basis of $\overline{V}$. Notice that the $K$-module $N$ is defined by $$N=\left\{g:\left(\bigcup\limits_{i\in \Delta}X_i\right)\cup I \rightarrow K: g(x)\neq 0 \text{ only for finitely many } x \right\}.$$ 
\fim

The main idea of theorem \ref{teorONequivrepres} was to give a condition over the graph $E$ that guarantees that the hypothesis of theorem \ref{theorem1} are satisfied. Below we give an example of a graph that does not satisfy the hypothesis of theorem \ref{theorem1}, yet its conclusion (and hence the conclusion of theorem \ref{teorONequivrepres}) is still valid.


\begin{exemplo} Consider the graph
 
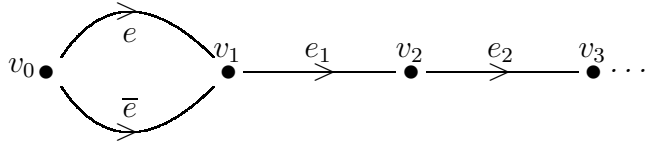
\centerline{
\setlength{\unitlength}{2cm}
\begin{picture}(4,0.6)
\put(0.0,0){\circle*{0.08}}
\put(-0.25,0.0){$v_0$}
\qbezier(0.1,-0.1)(0.5,-0.7)(1.1,-0.1)
\qbezier(0.1,0.1)(0.5,0.7)(1.1,0.1)
\put(0.45,-0.45){$>$}
\put(0.45,0.35){$>$}
\put(0.5,0.2){$e$}
\put(0.5,-0.3){$\overline{e}$}
\put(1.2,0){\circle*{0.08}}
\put(1.3,0){\line(1,0){1}}
\put(1.1,0.1){$v_1$}
\put(1.7,0.1){$e_1$}
\put(1.75,-0.05){$>$}
\put(2.4,0){\circle*{0.08}}
\put(2.5,0){\line(1,0){1}}
\put(2.3,0.1){$v_2$}
\put(2.9,0.1){$e_2$}
\put(2.9,-0.05){$>$}
\put(3.6,0){\circle*{0.08}}
\put(3.5,0.1){$v_3$}
\put(3.7,0){\dots}
\end{picture}}
\vspace{0.6cm}
\end{exemplo}
Notice that this graph does not satisfy the hypothesis of theorem \ref{theorem1} or the hypothesis of proposition \ref{prop1}. However, given a representation $\Phi:L_K(E)\rightarrow Hom_K(V)$, it is possible to choose basis of $V_e:=\Phi(e)\Phi(e^*)(V)$ and $V_u:=\Phi(u)(V)$ satisfying (\ref{h1}). Let us show how to choose such basis. First, fix a basis $B_{v_1}$ of $V_{v_1}$. Recall that for each $e\in E^0$, $\Phi(e):V_{r(e)}\rightarrow V_e$ and $\Phi(e^*):V_e\rightarrow V_{r(e)}$ are $K$-module isomorphisms. So, $B_e:=\Phi(e)(B_{v_1})$ is a basis of $V_e$ and $B_{\overline{e}}:=\Phi(\overline{e})(B_{v_1})$ is a basis of $V_{\overline{e}}$. Notice that $V_{e_i}=V_{v_i}$ for all $i\geq 1$. So, $B_{e_1}:=B_{v_1}$ is a basis of $V_{e_1}$. Define $B_{v_2}:=\Phi(e_1^*)(B_{e_1})$, which is a basis of $V_{v_2}$. Proceeding inductively, let $B_{e_i}:=B_{v_i}$, which is a basis of $V_{e_i}$ and define $B_{v_{i+1}}:=\Phi(e_{i}^*)(B_{e_i})$. This way we obtain basis satisfying condition (\ref{h1}).
Following theorem \ref{equivrep}, $\Phi$ is equivalent to a representation induced by an $E$-algebraic branching system.

\end{document}